# ON-LINE PREDICTIVE LINEAR REGRESSION[1]

By Vladimir Vovk, Ilia Nouretdinov and Alex Gammerman

*Royal Holloway, University of London*

We consider the on-line predictive version of the standard problem of linear regression; the goal is to predict each consecutive response given the corresponding explanatory variables and all the previous observations. The standard treatment of prediction in linear regression analysis has two drawbacks: (1) the classical prediction intervals guarantee that the probability of error is equal to the nominal significance level $\varepsilon$, but this property per se does not imply that the long-run frequency of error is close to $\varepsilon$; (2) it is not suitable for prediction of complex systems as it assumes that the number of observations exceeds the number of parameters. We state a general result showing that in the on-line protocol the frequency of error for the classical prediction intervals does equal the nominal significance level, up to statistical fluctuations. We also describe alternative regression models in which informative prediction intervals can be found before the number of observations exceeds the number of parameters. One of these models, which only assumes that the observations are independent and identically distributed, is popular in machine learning but greatly underused in the statistical theory of regression.

**1. Introduction.** Let $y_n$, $n = 1, 2, \ldots$, be the sequence of response variables to be predicted, and let $\mathbf{x}_n = (x_{n,1}, \ldots, x_{n,K})$, $n = 1, 2, \ldots$, be the corresponding vectors of explanatory variables. The standard assumption of linear regression analysis is that the explanatory vectors $\mathbf{x}_n$ are deterministic and

$$y_n = \alpha + \boldsymbol{\beta} \cdot \mathbf{x}_n + \xi_n, \quad (1)$$

where $\alpha$ is an unknown coefficient, $\boldsymbol{\beta} \in \mathbb{R}^K$ is an unknown vector of coefficients, and $\xi_n$, $n = 1, 2, \ldots$, are IID (independent and identically distributed)

Received October 2007; revised April 2008.
[1]Supported in part by EPSRC Grant EP/F002998/1, MRC Grant G0301107, the Cyprus Research Promotion Foundation and the Royal Society.

Gaussian random variables with mean 0 and unknown variance $\sigma^2 > 0$ [we will write $\xi_n \sim N(0, \sigma^2)$]. The model (1) will be called the *Gauss linear model*. It is the standard textbook model.

The standard classes of problems associated with the Gauss linear model are parameter estimation, testing hypotheses about parameters and prediction. In this paper we will be concerned only with prediction, mainly in the form of prediction intervals rather than point predictions.

A major drawback of the Gauss linear model is that the corresponding prediction intervals are uninformative (i.e., coincide with the whole real line) unless the number of observations exceeds the number of parameters. The responses of a complex system cannot be realistically expected to be modeled using a small number of parameters, whereas the number of observations can be very limited. This motivates consideration of three other models in this paper, none of which requires that the number of observations should exceed the number of parameters.

Perhaps the most important of these models is what we call the *IID model*: it is only assumed that the sequence of pairs $(\mathbf{x}_n, y_n)$ is IID. This model is nonparametric, effectively involving infinitely many parameters. Despite this, the model does allow one to obtain informative prediction intervals. The IID model, however, also has a fundamental limitation: informative prediction intervals become possible only when the number of observations reaches $1/\varepsilon$, where $\varepsilon$ is the chosen significance level.

Our third regression model combines the assumption (1) with the assumption that $\mathbf{x}_n$ are independent (between themselves and of $\xi_1, \xi_2, \ldots$) and identically distributed Gaussian random vectors. We call it the *MVA model*, with MVA referring to "multivariate analysis." It has also been widely discussed in the statistical literature; for example, Sampson's (1974) "two regressions" refers to the Gauss linear model and the MVA model. This model is narrower than both Gauss linear and IID models, and its strong assumptions ensure that informative prediction intervals can be produced almost right away.

Finally, we consider the combination of the Gauss linear and IID models, which we call the *IID–Gauss model*: in addition to (1) we assume that the explanatory vectors $\mathbf{x}_n$, $n = 1, 2, \ldots$, are random and IID (not necessarily Gaussian, as in the MVA model) and that the sequence $\xi_1, \xi_2, \ldots$ is independent of the explanatory vectors. This model, however, appears to be of secondary importance. Empirically, it allows informative prediction intervals at significance level $\varepsilon$ soon after the number of observations exceeds the minimum of $1/\varepsilon$ and the number of parameters.

All the models considered in this paper are shown in Figure 1, with arrows leading from more general to more specific models. In this paper we begin (in Section 5) with the IID model. This is the most common model used in modern day statistics and it does not involve the often unrealistic assumption that the noise variables $\xi_n$ are Gaussian or that the explanatory vectors $\mathbf{x}_n$ are Gaussian. An important advantage of the classical Gauss linear model, considered in Section 6, is that the explanatory vectors are not assumed to be IID (in other words, no "random design" is assumed). This model is essentially equivalent to making no assumptions whatsoever about the distribution of $\mathbf{x}_n$ and assuming that the $\xi_n$ in (1) are IID and distributed as $N(0, \sigma^2)$ conditional on $\mathbf{x}_1, \mathbf{x}_2, \ldots$. The Gauss linear model (understood in this way) and the IID model are not comparable between themselves, but both contain the other two models: the IID–Gauss model (Section 8), which is the intersection of the IID and Gauss linear models, and the MVA model (Section 7), which makes the further assumption that the explanatory vectors are Gaussian.

Fisher (1973), Section IV.3, emphatically defended the use of the Gauss linear model even in the case where the distribution of the explanatory vectors is known (with or without parameters). There is also a view in the literature that the Gauss linear model and the MVA model are "essentially equivalent" [for a review of some results in this direction, see Sampson (1974)]. Our conclusion, however, is similar to Brown's (1990): when the MVA model is true, it can be far more useful for prediction; in particular, it can start giving informative prediction intervals long before the number of observations reaches the number of parameters $K$ (or the inverse significance level $1/\varepsilon$).

This paper uses a general method of prediction called conformal prediction. The method is reviewed in detail in the monograph by Vovk, Gammerman and Shafer (2005) and introduced in the work leading up to that monograph. For each of the four models in Figure 1 we define a suitable confidence predictor, that is, a strategy for producing prediction intervals or, more generally, prediction regions. For the IID model we follow Vovk, Gammerman and Shafer (2005) and for the Gauss linear model we use Fisher's classical confidence predictor. The confidence predictors for the MVA and IID–Gauss models are new.

We are interested in two criteria of quality of confidence predictors, which we call "validity" and "accuracy." For valid confidence predictors, the prob-ability of error equals the nominal significance level $\varepsilon$ (or at least never exceeds $\varepsilon$, in which case we will refer to them as "conservatively valid," or just "conservative," confidence predictors). The second criterion is applied only to valid confidence predictors: we want the prediction intervals to be as narrow as possible; in this paper we, somewhat arbitrarily, measure the narrowness of a prediction interval $[a, b]$ by its length $b - a$. In particular, we want the prediction intervals to become bounded as soon as possible.

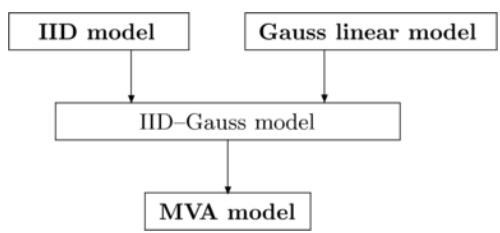

Fig. 1. *The four models considered in this paper (the three main models are given in boldface).*

Correspondingly, this paper uses two kinds of entities that one might want to call "models." The first kind is "hard models," such as the four models in Figure 1. These are the usual statistical models: our working hypothesis is that the data set was generated by one of the probability distributions in the model. In particular, the validity of our confidence predictors is allowed to depend on the hard model. By default, the word "model" means "hard model."

In addition to the accepted hard model, one often has other a priori information about the data-generating distribution: for example, only a few parameters might provide the bulk of the information relevant to prediction. Whereas we might hesitate to include such a priori information in the hard model explicitly, since it might destroy the validity of our confidence predictor if this information happened to be far from the truth, we might still be able to use such information in designing accurate confidence predictors provided our model is flexible enough. A running example in this paper, introduced in Section 4, will be a linear system with 100 parameters ten of which are felt to be especially important. This will be our "soft model" (not defined formally); whether it is true or not affects only the accuracy, but not validity, of our confidence predictors.

Separation of the available information about the data-generating distribution into the hard model and soft model increases robustness of confidence predictors with respect to modeling errors. If such an error occurs in the soft model, the validity of predictions is not affected. At worst the predictions will become useless, but they will not become misleading (with high probability under any distribution in the hard model). For a further discussion and empirical study, see Gammerman and Vovk (2007), Section 4.

The property of validity of conformal predictors can be stated in an especially strong form in the on-line prediction protocol. It turns out that the true responses fall outside the corresponding prediction regions independently for different observations. In combination with the law of large numbers this implies that, with high probability, the frequency of error is approximately equal to the nominal significance level. Surprisingly, even for the classical prediction intervals in the Gauss linear model this property had been unknown prior to the work leading up to Vovk, Gammerman and Shafer (2005).

Two recent reviews of the theory of conformal prediction are Gammerman and Vovk (2007) and Shafer and Vovk (2008). Parts of these papers are devoted to regression problems.

Section 2 formally introduces the on-line prediction protocol, with a more detailed discussion postponed until Section 9. In Section 3 we describe the method of conformal prediction and state two key results (proved in the Appendix): one asserts the strong validity and the other universality of conformal predictors. Section 4 describes an artificial data set used in later sections for illustrating the performance of various conformal predictors. The following Sections 5–8 apply the method of conformal prediction to the IID, Gauss linear, MVA and IID–Gauss models, in this order. Section 10 concludes.

**2. On-line protocol, part I.** In our prediction protocol, the task is to sequentially predict $y_n$, $n = 1, 2, \ldots$, from $\mathbf{x}_n$ and $(\mathbf{x}_i, y_i)$, $i = 1, \ldots, n-1$. This on-line protocol is popular in machine learning [see, e.g., Cesa-Bianchi and Lugosi (2006) and references therein], but most statistical research (except some work on sequential analysis) is still done in the "off-line," or "batch," framework, where one starts from a complete sample $(\mathbf{x}_1, y_1), \ldots, (\mathbf{x}_N, y_N)$. One of the few statisticians advocating the on-line protocol (under the name "prequential," or predictive sequential) has been Dawid (1984).

*Weak and strong validity and median accuracy.* To explain what precisely we mean by validity and accuracy, the two criteria of predictive performance mentioned in Section 1, we will need the notation introduced in the following description of the on-line prediction protocol.

ON-LINE PREDICTION PROTOCOL

FOR $n = 1, 2, \ldots$:
    Predictor observes $\mathbf{x}_n \in \mathbb{R}^K$;
    Predictor outputs $\Gamma_n^\varepsilon \subseteq \mathbb{R}$ for all $\varepsilon \in (0, 1)$;
    Predictor observes $y_n \in \mathbb{R}$;
    $\mathrm{err}_n^\varepsilon := \mathbb{I}_{y_n \notin \Gamma_n^\varepsilon}$ for all $\varepsilon \in (0, 1)$;
    $L_n^\varepsilon := \sup \Gamma_n^\varepsilon - \inf \Gamma_n^\varepsilon$ for all $\varepsilon \in (0, 1)$
END FOR.

(As usual, $\mathbb{I}_F$ is defined to be 1 if the condition $F$ holds and 0 if not.) At each step and for each significance level $\varepsilon$, Predictor outputs a *prediction region* (usually, although not necessarily, an interval) $\Gamma_n^\varepsilon \subseteq \mathbb{R}$. We require that, for all $n$, the family $\Gamma_n^\varepsilon$ of prediction regions should be nested: $\Gamma_n^{\varepsilon_1} \subseteq \Gamma_n^{\varepsilon_2}$ whenever $\varepsilon_1 > \varepsilon_2$. An error is registered, $\mathrm{err}_n^\varepsilon = 1$, if the prediction region fails to contain the true response $y_n$, and the accuracy of this particular prediction is measured by the length $L_n^\varepsilon$ of the corresponding *prediction interval* $\operatorname{co}\Gamma_n^\varepsilon$ ($\operatorname{co} E$ standing for the convex hull of the set $E$).

Let $\operatorname{Err}_n^\varepsilon := \operatorname{err}_1^\varepsilon + \cdots + \operatorname{err}_n^\varepsilon$ be the cumulative number of errors made up to, and including, step $n$. In the following sections, we will find it convenient to distinguish between two notions of validity, "weak validity" and "strong validity."

DEFINITION 1. A *confidence predictor* is defined to be a measurable prediction strategy $\Gamma_n^\varepsilon = \Gamma^\varepsilon(\mathbf{x}_1, y_1, \ldots, \mathbf{x}_{n-1}, y_{n-1}, \mathbf{x}_n)$ in the on-line prediction protocol.

DEFINITION 2. A confidence predictor is *weakly valid* in some statistical model if the probability that $\operatorname{err}_n^\varepsilon = 1$ is $\varepsilon$, for each $\varepsilon \in (0,1)$ and each $n$ under any probability distribution in the model.

The definition of weak validity is standard [cf. Cox and Hinkley (1974), (75) on page 243]. Weak validity by itself does not imply that $\operatorname{Err}_n/n$ is likely to be close to $\varepsilon$ for large $n$.

DEFINITION 3. A confidence predictor is *strongly valid* if it is weakly valid and, for each $\varepsilon \in (0,1)$, the events $\operatorname{err}_n^\varepsilon = 1$, $n = 1, 2, \ldots$, are independent.

Figure 3 below shows the plot of $\operatorname{Err}_n^\varepsilon$ against $n$ for a specific confidence predictor considered in this paper; it is typical of our predictors that the slopes of the plots of $\operatorname{Err}_n^\varepsilon$ are close to the corresponding significance levels $\varepsilon$ (we use the significance levels 5%, 1% and 0.5% in all our figures, represented by the corresponding *confidence levels* $1 - \varepsilon$ in the legends). This is the only figure in this paper illustrating the validity of our confidence predictors; such figures, in view of the mathematical results guaranteeing validity, tend to be uninformative.

We will measure the accuracy of the predictions made for the first $n$ observations by the median $M_n^\varepsilon$ of the sequence $L_1^\varepsilon, \ldots, L_n^\varepsilon$; again, this measure is arbitrary, to a large degree. A plot of $M_n^\varepsilon$ against $n$ will be called the *median-accuracy plot*; examples of such plots are given in Figures 2 and 4–6.

Unfortunately, the simple notions of validity introduced earlier have to be extended to become useful for our purpose. This is needed because, for example, the classical prediction intervals are uninformative before the number of observations reaches the number of parameters, and so for small $n$ the error probability is zero rather than $\varepsilon$. Let $\mathcal{N}$ be a set of positive integer numbers (we are mainly interested in the case where $\mathcal{N}$ has the form $\{m, m+1, \ldots\}$).

DEFINITION 4. We say that a confidence predictor is *weakly valid for* $n \in \mathcal{N}$ in a statistical model if the probability is $\varepsilon$ that it makes an error, $\text{err}_n^\varepsilon = 1$, at step $n$ under any probability distribution in the model and for all $n \in \mathcal{N}$ and $\varepsilon \in (0,1)$. It is *strongly valid for* $n \in \mathcal{N}$ if, in addition, $\text{err}_n^\varepsilon$, $n \in \mathcal{N}$, are independent for any fixed $\varepsilon$.

*The role of the on-line protocol.* The exposition of this paper is based on the on-line protocol, but the majority of our findings are not constrained to this specific protocol. For example, the fact that valid and informative prediction intervals can become feasible in the MVA model before the number of observations exceeds the number of parameters does not depend on the prediction protocol. In the absence of the on-line protocol, however, "validity" should be understood in the standard sense of weak validity.

**3. Conformal prediction.** In this section we define a class of confidence predictors, called conformal predictors, and state results about their validity and universality, in a certain sense.

*Notions of sufficiency.* Fix some *observation space* $Z$. We will be interested in the space $Z = \mathbb{R}^K \times \mathbb{R}$ of pairs $(\mathbf{x}, y)$; in general, $Z$ is a measurable space assumed to be Luzin, to ensure the existence of regular conditional probabilities. To define conformal predictors, we will need not only a statistical model on $Z^\infty$ but also a sequence of sufficient statistics $S_n : Z^n \to \Sigma_n$, $n = 1, 2, \ldots$; we will always assume that $\Sigma_n = S_n(Z^n)$. We will need a strengthened form of sufficiency; in our definitions we mainly follow Lauritzen (1988), Section II.2.

The sequence $(S_n)$ is *algebraically transitive* if there exists a sequence of measurable functions $F_n : \Sigma_{n-1} \times Z \to \Sigma_n$, $n = 2, 3, \ldots$, such that

$$S_n(\zeta_1, \ldots, \zeta_{n-1}, \zeta_n) = F_n(S_{n-1}(\zeta_1, \ldots, \zeta_{n-1}), \zeta_n)$$

for all $(\zeta_1, \ldots, \zeta_{n-1}, \zeta_n) \in Z^n$. Intuitively, $S_n(\zeta_1, \ldots, \zeta_n)$ is the summary of the first $n$ observations, and the condition of algebraic transitivity means that the summary can be updated on-line.

The sequence $(S_n)$ is *totally sufficient* for a statistical model $\mathcal{P}$ on $Z^\infty$ if, for each $n = 1, 2, \ldots$:

- $S_n$ is sufficient for $\mathcal{P}$;
- $\zeta_1, \ldots, \zeta_n$ and $\zeta_{n+1}, \zeta_{n+2}, \ldots$ are conditionally independent given $S_n(\zeta_1, \ldots, \zeta_n)$, where $(\zeta_1, \zeta_2, \ldots) \sim P$, for any $P \in \mathcal{P}$.

The second condition ensures that $S_n(\zeta_1, \ldots, \zeta_n)$ carries all information in $\zeta_1, \ldots, \zeta_n$ that can be used for predicting the future observations $\zeta_{n+1}, \zeta_{n+2}, \ldots$.

A sequence of statistics that is both algebraically transitive and totally sufficient will be called an *ATTS sequence*. In the rest of this paper we will often say "model" to mean a statistical model $\mathcal{P}$ equipped with an ATTS sequence $(S_n)$. This makes the word "model" ambiguous as we often omit "statistical" in "statistical model," but this should not lead to misunderstandings.

Each of the four statistical models considered in this paper (see Figure 1) will be complemented with an ATTS sequence; in all four cases the observation space $Z$ will be $\mathbb{R}^K \times \mathbb{R}$.

*Testing conformity.* The main ingredient of conformal prediction is statistical testing of conformity of a new observation $\zeta_n$ to the old observations $\zeta_1, \ldots, \zeta_{n-1}$. In general, our statistical tests will be randomized.

Fix a statistical model $\mathcal{P}$ with an ATTS sequence $S_n : Z^n \to \Sigma_n$. Define $\Sigma_0$ to be a fixed one-element set. Any sequence of measurable functions $A_n : \Sigma_{n-1} \times Z \to \mathbb{R}$, $n = 1, 2, \ldots$, is called a *nonconformity measure*; $A_n$ will be our test statistics. Given a nonconformity measure $(A_n)$, for each sequence $\zeta_1, \zeta_2, \ldots$ of observations and each sequence $\tau_1, \tau_2, \ldots \in [0,1]^\infty$ we define the *p-values*

$$\begin{aligned} p_n &= p_n(\zeta_1, \ldots, \zeta_n, \tau_n) \\ (2) \quad &:= \mathbb{P}(A_n^{\mathrm{rnd}} > A_n^{\mathrm{obs}} \mid S_n^{\mathrm{rnd}} = S_n^{\mathrm{obs}}) + \tau_n \mathbb{P}(A_n^{\mathrm{rnd}} = A_n^{\mathrm{obs}} \mid S_n^{\mathrm{rnd}} = S_n^{\mathrm{obs}}), \\ & \qquad\qquad\qquad\qquad\qquad\qquad\qquad\qquad\qquad\qquad n = 1, 2, \ldots, \end{aligned}$$

where $A_n^{\mathrm{rnd}} := A_n(S_{n-1}(\xi_1, \ldots, \xi_{n-1}), \xi_n)$ and $S_n^{\mathrm{rnd}} := S_n(\xi_1, \ldots, \xi_n)$ are the "random" values, $A_n^{\mathrm{obs}} := A_n(S_{n-1}(\zeta_1, \ldots, \zeta_{n-1}), \zeta_n)$ and $S_n^{\mathrm{obs}} := S_n(\zeta_1, \ldots, \zeta_n)$ are the "observed" values, and the probabilities are taken with respect to $(\xi_1, \xi_2, \ldots) \sim P$ for some $P \in \mathcal{P}$. Since $S_n$ are sufficient statistics, $p_n$ do not depend on $P \in \mathcal{P}$ (at least for a suitable choice of regular conditional probabilities). We will be interested in two cases: *deterministic*, where $\tau_n = 1$ for all $n$, and *randomized*, where $\tau_1, \tau_2, \ldots$ are generated independently from the uniform distribution $U$ on $[0,1]$ (such $\tau_1, \tau_2, \ldots$ model the output of a random numbers generator).

THEOREM 1. *Suppose that the sequence of observations $(\zeta_1, \zeta_2, \ldots) \in Z^\infty$ is generated from a probability distribution $P \in \mathcal{P}$ and that the random numbers $(\tau_1, \tau_2, \ldots) \sim U^\infty$ are independent of the observations. The p-values (2) are then independent and distributed uniformly on $[0,1]$:*

$$(p_1, p_2, \ldots) \sim U^\infty.$$

For a proof of this theorem, see the Appendix. The fact that $p_n \sim U$ is well known, at least in the continuous case [see, e.g., Cox and Hinkley (1974), page 66; (2) is a version of Cox and Hinkley's (1)].

*Conformal prediction.* We start by extending, and spelling out in a greater detail, the notion of a confidence predictor: in the general theory of this section and in its application to the IID model in Section 5 we will need an element (typically quite small) of randomization in confidence predictors.

DEFINITION 5. A *randomized confidence predictor* is a measurable function which maps every significance level $\varepsilon \in (0,1)$, every data sequence $\mathbf{x}_1, y_1, \ldots, \mathbf{x}_{n-1}, y_{n-1}$, every vector $\mathbf{x}_n$ of explanatory variables, and every number $\tau \in [0,1]$ to a set $\Gamma_n^\varepsilon = \Gamma^\varepsilon(\mathbf{x}_1, y_1, \ldots, \mathbf{x}_{n-1}, y_{n-1}, \mathbf{x}_n, \tau) \subseteq \mathbb{R}$. We will use the notation $\Gamma_n^\varepsilon$ when the data sequence, the vector of explanatory variables, and the number $\tau$ are clear from the context.

Let the observation space be $Z = \mathbb{R}^K \times \mathbb{R}$. Once the *p*-values (2) are defined, we can use them for confidence prediction [this is a standard procedure; cf. Cox and Hinkley (1974), (76) on page 243]: we set

$$\begin{aligned}(3)\quad &\Gamma^\varepsilon(\mathbf{x}_1, y_1, \ldots, \mathbf{x}_{n-1}, y_{n-1}, \mathbf{x}_n, \tau_n) \\ &:= \{y \in \mathbb{R} : p_n((\mathbf{x}_1, y_1), \ldots, (\mathbf{x}_{n-1}, y_{n-1}), (\mathbf{x}_n, y), \tau_n) > \varepsilon\}.\end{aligned}$$

DEFINITION 6. The randomized confidence predictor defined by (3) is called the *smoothed conformal predictor determined by* the nonconformity measure $(A_n)$. A *smoothed conformal predictor* is a smoothed conformal predictor determined by some nonconformity measure.

The following statement immediately follows from Theorem 1 and asserts that smoothed conformal predictors are strongly valid.

COROLLARY 1. *If the sequence of observations* $(\mathbf{x}_n, y_n)$, $n = 1, 2, \ldots$, *is generated by a probability distribution* $P \in \mathcal{P}$ *and a smoothed conformal predictor is fed with random numbers* $(\tau_1, \tau_2, \ldots) \sim U^\infty$ *independent of the observations, the error sequence* $\mathrm{err}_1^\varepsilon, \mathrm{err}_2^\varepsilon, \ldots$ *at any significance level $\varepsilon$ is a sequence of IID Bernoulli random variables with parameter $\varepsilon$.*

The adjective "smoothed" refers to using random numbers; if we take $\tau_n = 1$ for all $n = 1, 2, \ldots$, we will obtain the definition of a "deterministic conformal predictor," or just "conformal predictor," and in this case we omit $\tau_n$ from our notation.

DEFINITION 7. A *conformal predictor* is the confidence predictor defined by

$$\begin{aligned}&\Gamma^\varepsilon(\mathbf{x}_1, y_1, \ldots, \mathbf{x}_{n-1}, y_{n-1}, \mathbf{x}_n) \\ &:= \{y \in \mathbb{R} : p_n((\mathbf{x}_1, y_1), \ldots, (\mathbf{x}_{n-1}, y_{n-1}), (\mathbf{x}_n, y), 1) > \varepsilon\},\end{aligned}$$

where the *p*-values $p_n$ are defined by (2).

Notice that when a conformal predictor makes an error, the corresponding smoothed conformal predictor also makes an error. In combination with Corollary 1, we can see that conformal predictors are *conservative*, in the sense that, for each $\varepsilon$, their error sequence $\mathrm{err}_1^\varepsilon, \mathrm{err}_2^\varepsilon, \ldots$ is dominated by a sequence of IID Bernoulli random variables with parameter $\varepsilon$. In particular, whereas we have $\lim_{n\to\infty}(\mathrm{Err}_n^\varepsilon/n) = \varepsilon$ a.s. for smoothed conformal predictors, we only have $\limsup_{n\to\infty}(\mathrm{Err}_n^\varepsilon/n) \le \varepsilon$ a.s. for conformal predictors.

We will see that there is no difference between conformal predictors and the corresponding smoothed conformal predictors for the Gauss linear model and $n \ge K + 3$ since the second addend on the right-hand side of (2) is then zero. There is also no difference for the MVA model and $n \ge 3$; however, the difference is important (although usually barely noticeable on error and accuracy plots) for the IID model.

A natural question is whether there are other ways to achieve validity, except conformal prediction. The following theorem will give a negative answer to a version of this question.

DEFINITION 8. A confidence predictor $\Gamma$ is *invariant* if $\Gamma_n^\varepsilon$, $n > 1$, depends on the first $n-1$ observations only through the value of $S_{n-1}$ on those observations.

The use of invariant confidence predictors is natural in view of the sufficiency principle; see, for example, Cox and Hinkley (1974), Section 2.3(ii). Let $\mathcal{N}$ be a set of positive integers. We say that a confidence predictor $\Gamma^\dagger$ is *at least as accurate as* another confidence predictor $\Gamma$ for $n \in \mathcal{N}$ if

$$(\Gamma^\dagger)^\varepsilon(\mathbf{x}_1, y_1, \ldots, \mathbf{x}_{n-1}, y_{n-1}, \mathbf{x}_n) \subseteq \Gamma^\varepsilon(\mathbf{x}_1, y_1, \ldots, \mathbf{x}_{n-1}, y_{n-1}, \mathbf{x}_n)$$

for all $\varepsilon$, all $n \in \mathcal{N}$, and $P$-almost all $\mathbf{x}_1, y_1, \ldots, \mathbf{x}_{n-1}, y_{n-1}, \mathbf{x}_n$, under any probability distribution $P \in \mathcal{P}$.

Recall that a statistic $S$ taking values in a measurable space $\Sigma$ is said to be *boundedly complete* (with respect to the statistical model $\mathcal{P}$) if, for any bounded measurable function $f: \Sigma \to \mathbb{R}$, the following condition is satisfied: the expected value $\mathbb{E}_P(f(S))$ of $f(S)$ is zero under all $P \in \mathcal{P}$ only if $f(S) = 0$ $P$-almost surely for all $P \in \mathcal{P}$.

THEOREM 2. *Let $\mathcal{N}$ be a set of positive integers. Suppose the ATTS statistics $S_n$ are boundedly complete for $n \in \mathcal{N}$. If a confidence predictor $\Gamma$ is invariant and weakly valid for $n \in \mathcal{N}$, then there is a conformal predictor that is at least as accurate as $\Gamma$ for $n \in \mathcal{N}$.*

This theorem is also proved in the Appendix. An important step toward its proof was made by Takeuchi (1975), page 31.

TABLE 1
*Steps at which informative prediction becomes possible for the four models; $\varepsilon$ is the significance level ($\varepsilon < 1/2$ is assumed) and $K$ is the number of parameters*

| Model | The first step at which prediction intervals can become informative |
|---|---|
| IID model | $\lceil 1/\varepsilon \rceil$ |
| Gauss linear model | $K + 3$ |
| MVA model | 3 |
| IID–Gauss model | $\min(\lceil 1/\varepsilon \rceil, K + 3)$ |

The condition of bounded completeness holds for the Gauss linear model and the MVA model by the standard completeness result for exponential statistical models [see, e.g., Theorem 4.1 in Lehmann (1986)], and it is also known to hold for the IID model [see the theorem on page 797 in Bell, Blackwell and Breiman (1960)].

**4. Data set.** We will illustrate the accuracy of various confidence predictors using the following artificially generated data set with 600 observations and $K = 100$ explanatory variables. The components $x_{n,k}$ of $\mathbf{x}_n$ are independently generated from $N(0,1)$, and the responses $y_n$ are generated according to (1) with $\xi_n \sim N(0,1)$ independent between themselves and of all $x_{n,k}$, with $\alpha = 100$ and with the following components $\beta_k$ of $\boldsymbol{\beta}$:

$$\beta_k := \begin{cases} (-1)^{k-1} 10, & k = 1, \ldots, 10, \\ (-1)^{k-1}, & k = 11, \ldots, 100. \end{cases}$$

The probability distribution generating this data set belongs to all four models considered in this paper (Figure 1). It is natural to expect that more specific models, when true, will lead to better predictions. In one respect this is true: more general models allow informative predictions later, as shown in Table 1 (to be explained in later sections). However, soon after the threshold given in the table is reached, the quality of prediction becomes very similar on our data set.

The (informal) soft model guiding the choice of the nonconformity measure will include the assumption of linearity (1) and the knowledge, or guess, that the first 10 explanatory variables are much more important than the rest.

Relationship (1) between the response and explanatory variables can be written as

(4) $$y_n = \boldsymbol{\gamma} \cdot \mathbf{z}_n + \xi_n,$$

where

$$\boldsymbol{\gamma} := \begin{pmatrix} \alpha \\ \boldsymbol{\beta} \end{pmatrix} \in \mathbb{R}^{K+1} \quad \text{and} \quad \mathbf{z}_n := \begin{pmatrix} 1 \\ \mathbf{x}_n \end{pmatrix} \in \mathbb{R}^{K+1}.$$

For $l = 1, 2, \ldots$, let $\mathbf{Z}_l$ be the $l \times (K+1)$ matrix whose rows are $\mathbf{z}'_i$, $i = 1, \ldots, l$, and $\mathbf{y}_l$ be the vector whose $i$th element is $y_i$, $i = 1, \ldots, l$. We will sometimes refer to the first column of $\mathbf{Z}_l$ as the *dummy* column.

**5. The IID model.** The statistical model considered in this section is nonparametric: we simply assume that the observations $(\mathbf{x}_n, y_n)$ are IID. Notice that this does not involve the assumption of linearity of the "true" regression function or the assumption of a Gaussian noise. Linearity is, however, an important component of the soft model used for choosing a suitable nonconformity measure.

The ATTS statistics are

$$S_n := \{(\mathbf{x}_1, y_1), \ldots, (\mathbf{x}_n, y_n)\},$$

where we use $\{a_1, \ldots, a_n\}$ to denote the bag, or multiset, consisting of $a_1, \ldots, a_n$ (some of these elements may coincide). For each $n$, the conditional distribution of $(\xi_1, \ldots, \xi_n)$ given that

$$\{\xi_1, \ldots, \xi_n\} = \{(\mathbf{x}_1, y_1), \ldots, (\mathbf{x}_n, y_n)\},$$

where $\xi_i$ are IID random elements taking values in $\mathbb{R}^K \times \mathbb{R}$, assigns (with probability one) the same probability, $1/n!$, to every ordering $(\mathbf{x}_{\pi(1)}, y_{\pi(1)}), \ldots, (\mathbf{x}_{\pi(n)}, y_{\pi(n)})$ of the bag $\{(\mathbf{x}_1, y_1), \ldots, (\mathbf{x}_n, y_n)\}$.

The IID model is typical in that there is a great flexibility in choosing a nonconformity measure for use in conformal prediction. Suppose, for example, that the number of explanatory variables $K$ is too large for us to estimate all the $\beta_k$ and $\alpha$ in the soft model (1). We believe, however, that the first $K_n^\dagger \ll K$ of the explanatory variables are especially important, and it is feasible to estimate the corresponding $\beta_k$, $k = 1, \ldots, K_n^\dagger$, and $\alpha$.

Fix temporarily a positive integer number $n$. We will write $\mathbf{y}$ for $\mathbf{y}_n$, $\mathbf{Z}$ for $\mathbf{Z}_n$ and $K^\dagger$ for $K_n^\dagger$. Let $\mathbf{U}$ be the submatrix of $\mathbf{Z}$ consisting of the first $K^\dagger + 1$ columns of $\mathbf{Z}$: those that correspond to the explanatory variables deemed to be useful at this stage plus the dummy column $\mathbf{1}$. To test the conformity of the $n$th observation to the first $n - 1$ observations, we will first fit a hyperplane to all $n$ observations using the relevant explanatory variables. Applying a small "ridge coefficient" $a > 0$ to avoid the need to invert singular matrices, we obtain the vector of residuals

(5)  $$\mathbf{e} := \mathbf{y} - \mathbf{U}(\mathbf{U}'\mathbf{U} + a\mathbf{I})^{-1}\mathbf{U}'\mathbf{y},$$

whose components will be denoted $e_1, \ldots, e_n$.

We will be interested in the conformal predictor determined by the nonconformity measure

(6) $$A_n(S_{n-1}(\mathbf{x}_1, y_1, \ldots, \mathbf{x}_{n-1}, y_{n-1}), (\mathbf{x}_n, y_n)) := |e_n|.$$

Deleted and, especially, studentized residuals would also be a natural choice [see, e.g., Vovk, Gammerman and Shafer (2005), pages 34–35]. In our experience, however, the difference is not significant, and we stick to the simplest choice. The confidence predictor obtained from this conformal predictor by replacing the prediction regions $\Gamma_n^\varepsilon$ with the prediction intervals $\operatorname{co}\Gamma_n^\varepsilon$ will be called the *IID predictor* (cf. the comments at the end of this section).

The IID predictor can be implemented fairly efficiently. First notice that for the IID model the formula (2) for $p$-values can be simplified to

(7) $$p_n = \frac{|\{i : \alpha_i > \alpha_n\}| + \tau_n |\{i : \alpha_i = \alpha_n\}|}{n},$$

where $\alpha_i := A_n(\langle \zeta_1, \ldots, \zeta_{i-1}, \zeta_{i+1}, \ldots, \zeta_n \rangle, \zeta_i)$, $i$ ranges over $\{1, \ldots, n\}$, and $|E|$ stands for the size of the set $E$. In the case of the nonconformity measure (6), $\alpha_i = |e_i|$. The residuals (5) can be written in the form

$$\mathbf{e} = \mathbf{y} - \mathbf{U}(\mathbf{U}'\mathbf{U} + a\mathbf{I})^{-1}\mathbf{U}'\mathbf{y} = \mathbf{C}\mathbf{y},$$

where $\mathbf{C}$ is the matrix $\mathbf{I} - \mathbf{U}(\mathbf{U}'\mathbf{U} + a\mathbf{I})^{-1}\mathbf{U}'$, not depending on the response variables. If we fix the first $n-1$ response variables $y_i$ and vary the last one, $y$, the residuals $e_i = e_i(y)$, $i = 1, \ldots, n$, become linear functions of $y$ (this fact will also be used in Section 7). By (7) with $\tau_n := 1$, the $p$-value is the fraction of $i = 1, \ldots, n$ satisfying $|e_i(y)| \geq |e_n(y)|$; therefore, as $y$ varies from $-\infty$ to $\infty$, the $p$-value can change only at the at most $2n - 2$ points (called *critical* points) which are solutions to the linear equations $e_i(y) = e_n(y)$ and $e_i(y) = -e_n(y)$. This divides the real line into at most $4n - 3$ intervals: the critical points, considered as degenerate closed intervals, the open intervals bounded on both sides by adjacent critical points, and the two unbounded open intervals to the left of the leftmost critical point and to the right of the rightmost critical point; if there are no critical points, this collapses into one unbounded open interval $\mathbb{R}$. We can compute the $p$-value for one point in each of these intervals and then compute $\Gamma_n^\varepsilon$ as the union of the intervals with $p$-values exceeding $\varepsilon$. The computation of the IID prediction interval $\operatorname{co}\Gamma_n^\varepsilon$ can be simplified if we notice that the set $\Gamma_n^\varepsilon$ is closed (which is opposite to what we will have for the Gauss linear and MVA models): assuming that the set of critical points is nonempty, $\operatorname{co}\Gamma_n^\varepsilon$ is bounded if and only if the two unbounded intervals have $p$-values at most $\varepsilon$, in which case the end-points of $\operatorname{co}\Gamma_n^\varepsilon$ can be found as the leftmost and rightmost critical points with $p$-values exceeding $\varepsilon$. Computing $\Gamma_n^\varepsilon$ and $\operatorname{co}\Gamma_n^\varepsilon$ from scratch (e.g., without using the results of computations from the previous steps of the on-line protocol) takes time $O(n \log n)$ [see Vovk, Gammerman and Shafer (2005), page 33].

For use in our experiments with the artificial data set described in Section 4, we take

(8) $$K_n^\dagger := \begin{cases} 10, & \text{if } n < 103, \\ 100, & \text{otherwise,} \end{cases}$$

and so define $\mathbf{U}$ as the first 11 columns of $\mathbf{Z}$ if $n < 103$ and as the full $\mathbf{Z}$ otherwise. Our chosen value for the threshold, 103, appeared to us slightly less arbitrary than other choices, since it is the first step when the classical prediction intervals [see (10)] become bounded. However, the quality of the estimates of $\alpha$ and the 100 components of $\boldsymbol{\beta}$ is still poor when $n$ is close to 103. This affects the quality of our prediction intervals but does not show on the median-accuracy plots. The value of the ridge coefficient is always $a = 0.01$.

As Figure 2 shows, the IID predictor works well for our data set if the significance level is not too demanding: it can be seen from (7) (with $\tau_n := 1$) that for the IID prediction interval $\operatorname{co} \Gamma_n^\varepsilon$ to be bounded the number of observations $n$ has to be at least $1/\varepsilon$ (as Table 1 says). For example, for the significance level $\varepsilon = 0.5\%$, the IID predictor requires 200 observations to produce bounded predictions, and this shows on the median-accuracy plot at $n = 399$ (since for $n < 399$ at least half of the observed prediction intervals are infinitely wide).

The IID model is nonparametric but we can see that it still admits valid confidence predictors (or conservative confidence predictors if one insists on using deterministic predictors). The threshold $1/\varepsilon$ can be said to play the role of the number of parameters, and the nonparametric nature of the model is reflected in the fact that $1/\varepsilon \to \infty$ as $\varepsilon \to 0$. Since $1/\varepsilon$ tends to $\infty$ relatively slowly, such an infinite-dimensional model may be better for the purpose of prediction than a $K$-dimensional model with a very large $K$.

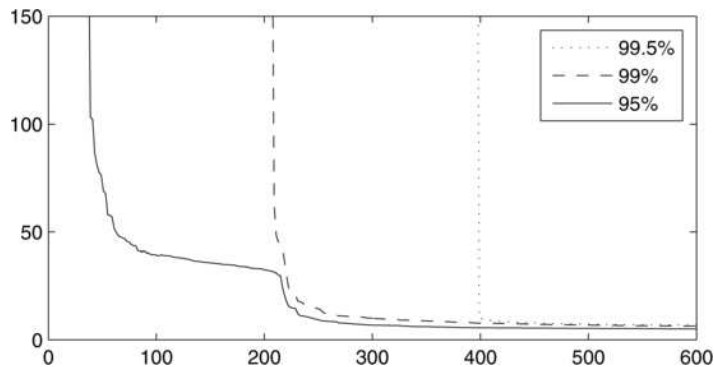

Fig. 2. *The median-accuracy plot for the IID predictor. The three significance levels used in this and all the following figures are $\varepsilon = 0.05, 0.01, 0.005$, shown in the form $100(1-\varepsilon)\%$ (the corresponding confidence levels) in the legends.*

ON-LINE REGRESSION                     15

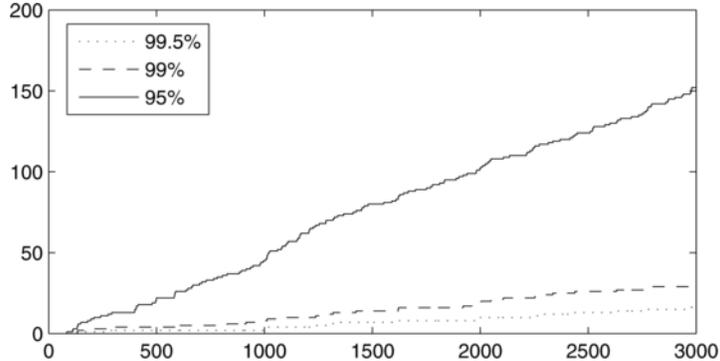

FIG. 3. *The cumulative numbers of errors made by the IID predictor:* $\mathrm{Err}_n^\varepsilon$ *is plotted against* $n$.

Theorem 2 is not directly applicable to the IID model, since only smoothed conformal predictors are valid, as the latter term is used in this paper. Vovk, Gammerman and Shafer (2005), Section 2.4, state two results of the same nature about the IID model.

There are two sources of conservativeness for the IID predictor as described above (and used for producing Figure 2). First, we used a deterministic predictor (taking $\tau_n = 1$ for all $n$), and second, we replaced each prediction region by its convex hull. Our experiments (see, e.g., Figure 3) show that we still have approximate validity.

For each model considered in this paper except the Gauss linear model we define a nonconformity measure involving the matrix $\mathbf{U}$ defined earlier in this section. In the case of the IID model, we have used the nonconformity measure (6) and called the corresponding conformal predictor with $\Gamma_n^\varepsilon$ replaced by $\mathrm{co}\,\Gamma_n^\varepsilon$ the IID predictor [it was called "Ridge Regression Confidence Machine" in Vovk, Gammerman and Shafer (2005)]. Of course, our brief term is somewhat misleading: it should always be borne in mind that the conformal predictor leading to the IID predictor is only one of many conformal predictors that can be defined in the IID model. Similarly, in the following three sections we will introduce the Gauss predictor, the MVA predictor and the IID–Gauss predictor, which will also correspond to specific nonconformity measures.

**6. The Gauss linear model.** Let $\hat{\boldsymbol{\gamma}}_l := (\mathbf{Z}_l'\mathbf{Z}_l)^{-1}\mathbf{Z}_l'\mathbf{y}_l$ be the least-squares estimate of the parameter vector $\boldsymbol{\gamma}$ in (4) from the first $l$ observations. For simplicity, we will assume that the matrix $\mathbf{Z}_l$ has full rank [i.e., $\mathrm{rank}\,\mathbf{Z}_l = \min(l, K+1)$] for all $l$; this implies that $\hat{\boldsymbol{\gamma}}_l$ is well defined for $l \geq K + 1$.

Let $\hat{y}_n$ be the least-squares prediction $\hat{\boldsymbol{\gamma}}_{n-1} \cdot \mathbf{z}_n$ for $y_n$ and

$$\hat{\sigma}_l^2 := \frac{1}{l - K - 1}(\mathbf{y}_l - \mathbf{Z}_l\hat{\boldsymbol{\gamma}}_l)'(\mathbf{y}_l - \mathbf{Z}_l\hat{\boldsymbol{\gamma}}_l)$$

be the standard estimate of $\sigma^2$ from $\mathbf{Z}_l$ and $\mathbf{y}_l$. It is well known that in the Gauss linear model the ratio

$$(9) \quad T_n := \frac{y_n - \hat{y}_n}{\sqrt{1 + \mathbf{z}_n'(\mathbf{Z}_{n-1}'\mathbf{Z}_{n-1})^{-1}\mathbf{z}_n}\hat{\sigma}_{n-1}}, \qquad n = K+3, K+4, \ldots,$$

has the $t$-distribution with $n - K - 2$ degrees of freedom. This gives the classical weakly valid prediction interval for the $n$th response,

$$(10) \quad \Gamma_n^\varepsilon := \{y \in \mathbb{R} : |y - \hat{y}_n| < t_{n-K-2}^{\varepsilon/2}\sqrt{1 + \mathbf{z}_n'(\mathbf{Z}_{n-1}'\mathbf{Z}_{n-1})^{-1}\mathbf{z}_n}\hat{\sigma}_{n-1}\},$$

$$n \geq K + 3,$$

where $t_m^\delta$ is the upper $\delta$ point of the $t$-distribution with $m$ degrees of freedom. [See, e.g., Seber and Lee (2003), (5.27).] We set $\Gamma_n^\varepsilon$ to $\mathbb{R}$ when $n < K + 3$.

Later in this section we will see that Corollary 1 implies the following property of the classical prediction intervals for the Gauss linear model.

COROLLARY 2. *Let $\varepsilon \in (0,1)$. The events $y_n \notin \Gamma_n^\varepsilon$, $n = K+3, K+4, \ldots$, are independent. In particular, the confidence predictor (10) is strongly valid for $n \geq K + 3$.*

REMARK. We have not seen Corollary 2 stated explicitly in the literature, but some closely related facts are known. Lemma 1 in Brown, Durbin and Evans (1975) asserts that (9) with $\hat{\sigma}_{n-1}$ removed are independent $N(0, \sigma^2)$ random variables; this can be used for prediction when the standard deviation $\sigma$ is known. Seillier-Moiseiwitsch [(1993), Example 1] shows that the statistics $T_n$ are independent when $K = 0$. It is interesting that both papers use the independence of $T_n$ for testing rather than for prediction.

Let us now see that some conformal predictor outputs the classical prediction intervals (10). This will demonstrate that Corollary 2 is indeed a special case of Corollary 1.

The ATTS statistics for the Gauss linear model are

$$S_n(\mathbf{x}_1, y_1, \ldots, \mathbf{x}_n, y_n) := \left(\mathbf{x}_1, \ldots, \mathbf{x}_n, \sum_{i=1}^n y_i, \sum_{i=1}^n y_i \mathbf{x}_i, \sum_{i=1}^n y_i^2\right).$$

(It is natural to have $\mathbf{x}_1, \ldots, \mathbf{x}_n$ as components of $S_n$, although they are superfluous under our original definition, in which $\mathbf{x}_1, \mathbf{x}_2, \ldots$ are deterministic.) The prediction intervals (10) are precisely the prediction regions output by the conformal predictor corresponding to the nonconformity measure

$$(11) \quad A_n(S_{n-1}(\mathbf{x}_1, y_1, \ldots, \mathbf{x}_{n-1}, y_{n-1}), (\mathbf{x}_n, y_n))$$

$$:= \frac{|y_n - \hat{y}_n|}{\sqrt{1 + \mathbf{z}_n'(\mathbf{Z}_{n-1}'\mathbf{Z}_{n-1})^{-1}\mathbf{z}_n}\hat{\sigma}_{n-1}}$$

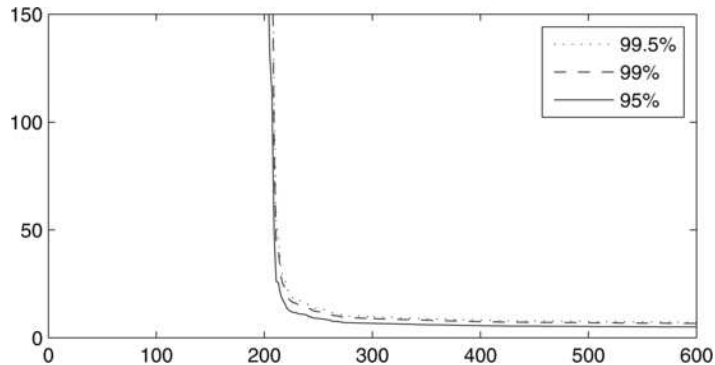

Fig. 4. *The median-accuracy plot for the classical prediction intervals.*

[cf. (9); the goodness of the definition follows from the formulas given at the beginning of this section]. The expression on the right-hand side of (11) can be replaced by other natural expressions, such as $|y_n - \hat{y}_n|$. See Vovk, Gammerman and Shafer (2005), Section 8.5, for further details.

According to our general convention, the conformal predictor (10) is called the *Gauss predictor* (although its discoverer was Fisher rather than Gauss).

We have already mentioned that the classical confidence predictor, $\Gamma_n^\varepsilon$ given by (10), does not work when there are many parameters; in particular, it is required that $n \geq K + 3$. Theorem 2 shows that there is hardly any way to use the knowledge that the first 10 explanatory variables are the important ones without abandoning the Gauss linear model: no weakly valid confidence predictor in a very wide and natural class can produce informative prediction intervals unless $n \geq K + 3$. Indeed, since the conditional distribution of the first $n$ observations given $S_n$ is concentrated at one point for $n \leq K + 1$ and at two points for $n = K + 2$ with probability one, no conformal predictor and, therefore, no weakly valid invariant confidence predictor can give a bounded prediction region $\Gamma_n^\varepsilon$ for $\varepsilon < 0.5$ and $n \leq K + 2$.

REMARK. A common reaction to the importance of the condition $n \geq K + 3$ is that one can use only a subset of explanatory variables when $n < K + 3$. We are, however, interested in confidence predictors that are valid under the Gauss linear model (1), not under some other model that is only "approximately true," in some ill-defined sense.

Figure 4 gives the median-accuracy plot for the confidence predictor (10); the predictor works very well soon after the number of observations reaches $K + 3 = 103$. Since the median is plotted, the good quality of the prediction intervals shows only from $n = 205$: indeed, for $n < 205$ at least half of the observed prediction intervals are infinitely wide.

**7. The MVA model.** Remember that the MVA model assumes, besides (1), that $\mathbf{x}_n$ are generated independently from the same unknown multivariate Gaussian distribution on $\mathbb{R}^K$, with the noise random variables $\xi_1, \xi_2, \ldots$ independent of $\mathbf{x}_1, \mathbf{x}_2, \ldots$. The ATTS statistics in the MVA model are

$$S_n := \left( \sum_{i=1}^n \mathbf{x}_i, \sum_{i=1}^n y_i, \sum_{i=1}^n \mathbf{x}_i \mathbf{x}_i', \sum_{i=1}^n y_i \mathbf{x}_i, \sum_{i=1}^n y_i^2 \right);$$

equivalently, the ATTS statistics can be defined to be the empirical means and covariances of all variables, that is, the response and the explanatory variables.

Let $\mathbf{y} := \mathbf{y}_n$, $\mathbf{Z} := \mathbf{Z}_n$, $K^\dagger := K_n^\dagger$ and $\mathbf{U}$ be as in Section 5. Suppose the value of the statistic $S_n$ is known. The vector of residuals (5) can now be written as

$$(12) \quad \mathbf{e} := \mathbf{y} - \mathbf{U}(\mathbf{U}'\mathbf{U} + a\mathbf{I})^{-1}\mathbf{U}'\mathbf{y} = \mathbf{y} - \mathbf{U}\mathbf{c},$$

where $\mathbf{c} := (\mathbf{U}'\mathbf{U} + a\mathbf{I})^{-1}\mathbf{U}'\mathbf{y}$ is a known vector. Since the joint distribution of $\mathbf{y}$ and the nondummy columns of $\mathbf{U}$ is invariant with respect to rotations around the vector $\mathbf{1}$, the distribution of $\mathbf{e}$ will also be invariant with respect to such rotations. It might help the reader's intuition to notice that knowing the value of $S_n$ is equivalent to knowing the lengths of and the angles between the following $K+2$ vectors: the $K+1$ columns of $\mathbf{Z}$ and $\mathbf{y}$.

In the rest of this section we will assume $n \geq 3$ (with arbitrary conventions for $n = 1, 2$). Let $e_1, \ldots, e_n$ be the components of the vector (12) of residuals and $\overline{e}_{n-1}$ be the average of $e_1, \ldots, e_{n-1}$. A standard statistical result [Fisher (1925)] allows us to conclude that

$$(13) \quad \sqrt{\frac{n-1}{n}} \frac{e_n - \overline{e}_{n-1}}{\sqrt{(1/(n-2)) \sum_{i=1}^{n-1} (e_i - \overline{e}_{n-1})^2}}$$

has the $t$-distribution with $n - 2$ degrees of freedom.

Let us see how to implement the conformal predictor corresponding to the nonconformity measure

$$(14) \quad A_n(S_{n-1}(\mathbf{x}_1, y_1, \ldots, \mathbf{x}_{n-1}, y_{n-1}), (\mathbf{x}_n, y_n)) := \frac{e_n - \overline{e}_{n-1}}{\sqrt{\sum_{i=1}^{n-1} (e_i - \overline{e}_{n-1})^2}},$$

which is proportional to (13); the fact that the right-hand side of (14) depends on the first $n - 1$ observations only through the value of $S_{n-1}$ can be seen from the representation (12), where $\mathbf{c}$ is a known vector. First we replace the true value $y_n$ by variable $y$ ranging over $\mathbb{R}$. Each residual $e_i$ becomes a linear [according to (12), where $\mathbf{c}$ also depends on $y$] function $e_i(y)$ of $y$, and the prediction region can be written as

$$\Gamma_n^\varepsilon := \left\{ y \in \mathbb{R} : \sqrt{\frac{n-1}{n}} \frac{|e_n(y) - \overline{e}_{n-1}(y)|}{\sqrt{(1/(n-2)) \sum_{i=1}^{n-1} (e_i(y) - \overline{e}_{n-1}(y))^2}} < t_{n-2}^{\varepsilon/2} \right\}.$$

The inequality in this formula is quadratic in $y$, so $\Gamma_n^\varepsilon$ is easy to find. We can see that the prediction region for $y_n$ is an interval (empirically, this is the typical case), the union of two rays, the empty set, or the whole real line.

Replacing $\Gamma_n^\varepsilon$ by $\operatorname{co}\Gamma_n^\varepsilon$ in the conformal predictor we have just defined gives the *MVA predictor*. Our experiments with the artificial data set of Section 4 are carried out as before [cf. (8)]: $\mathbf{U}$ is defined as the first 11 columns of $\mathbf{Z}$ if $n < 103$ and as the full $\mathbf{Z}$ otherwise.

The median-accuracy plot for the MVA predictor and our artificial data set is shown in Figure 5. Before the threshold 103 the predictor quickly learns $\alpha$ and the first 10 parameters $\beta_k$, and its performance more or less stabilizes before quickly improving again when it starts learning the other parameters from $n = 103$ onward; the second improvement in the performance shows on the median-accuracy plot from $n = 205$.

The performance of the MVA predictor is better than the performance of any other confidence predictor considered in this paper. Of course, this should not be taken to mean that the other predictors are worse. Different predictors are based on different information about the data set. None of the predictors "knows" that the components of $\mathbf{x}_n$ are realizations of independent standard Gaussian random variables; even the MVA model, the narrowest model considered in this paper, allows arbitrary means of and arbitrary correlations between different explanatory variables for the same observation. The Gauss predictor does not know that the $\mathbf{x}_n$ are IID and Gaussian. The IID predictor only knows that the observations $(\mathbf{x}_n, y_n)$ are IID, and the IID–Gauss predictor, introduced in the next section, knows, in addition, that the $y_n$ are generated by (1).

The median-accuracy plot for each of the four predictors is essentially determined by that for the MVA predictor and the threshold for the corresponding model as shown in Table 1. It is convenient to represent each line on a median-accuracy plot as the function that maps each value for the accuracy in the interval $[0, 150]$ to the first step at which that accuracy is achieved (so the graph of this function is obtained by rotating the page by $90°$ counterclockwise). Each of the three functions in Figure 2 is, approximately, the maximum of $2\lceil 1/\varepsilon \rceil$ and the corresponding function in Figure 5. Similarly, each of the three functions in Figure 4 is, approximately, the maximum of $2(K+3) = 206$ and the corresponding function in Figure 5. As usual, the factor of 2 appears because of the use of median in our accuracy plots.

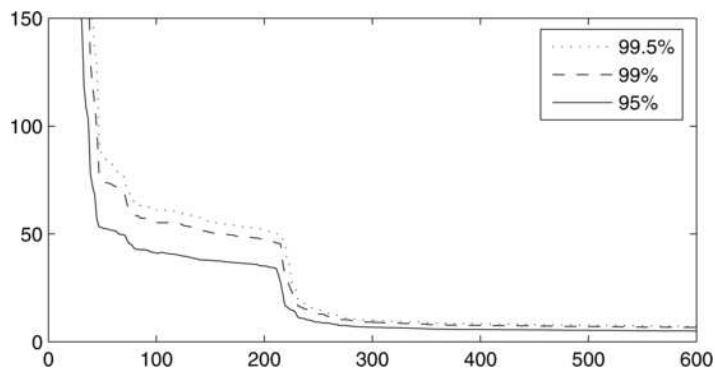

FIG. 5. *The median-accuracy plot for the MVA predictor.*

**8. The IID–Gauss model.** As defined in Section 1, the IID–Gauss model is the combination of the Gauss linear and IID models: we assume both that the observations are IID and that the responses are generated by (1) with $\xi_1, \xi_2, \ldots$ independent of $\mathbf{x}_1, \mathbf{x}_2, \ldots$. Correspondingly, the ATTS statistics are

$$S_n := \left( \{\mathbf{x}_1, \ldots, \mathbf{x}_n\}, \sum_{i=1}^n y_i, \sum_{i=1}^n y_i \mathbf{x}_i, \sum_{i=1}^n y_i^2 \right).$$

Using the nonconformity measure (6) and replacing the prediction regions output by the corresponding conformal predictor with their convex hulls, we obtain the *IID–Gauss predictor*. Its performance on our usual data set is shown in Figure 6. We do not know whether the IID–Gauss predictor can be implemented efficiently, and Figure 6 was produced using Monte-Carlo sampling from the conditional distributions given $S_n$. However, comparing Figure 6 to Figures 2 (to the left of $n = 205$) and 4 (to the right of $n = 205$), we can see that the following simple confidence predictor will work almost as well as the IID–Gauss predictor on our data set: predict using the IID predictor if $n < 103$ and predict using the Gauss predictor if $n \geq 103$. As in all other cases in this paper where the threshold $n = K + 3 = 103$ appears, the best switch-over point will be slightly greater than $K+3$, but the question of when exactly to switch is outside the scope of this paper.

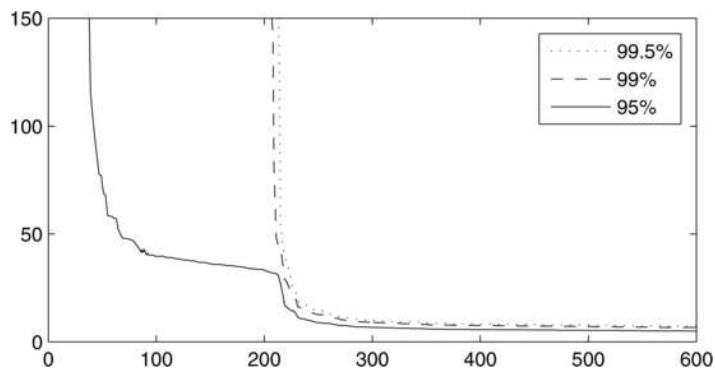

Fig. 6. *The median-accuracy plot for the IID–Gauss predictor.*

REMARK. The IID predictor and the IID–Gauss predictor use the same nonconformity measure, (6), but still produce very different median-accuracy plots at confidence level 99.5%. This happens because of the conditioning on the event $S_n^{\text{rnd}} = S_n^{\text{obs}}$ in the definition (2). Since the ATTS statistics perform more radical data compression in the case of the IID–Gauss model, the achievable values of $\mathbb{P}(A_n^{\text{rnd}} \geq A_n^{\text{obs}} \mid S_n^{\text{rnd}} = S_n^{\text{obs}})$ [corresponding to (2) with $\tau_n := 1$] are much smaller than the $1/n$ achievable under the IID model.

As in the previous section, there is a close connection between Figures 5 and 6: each of the three functions in Figure 6 is, approximately, the maximum of $2\min(\lceil 1/\varepsilon \rceil, K+3)$ and the corresponding function in Figure 5. The distributive law of max over min now implies that each of the three functions in Figure 6 is the minimum of the corresponding functions in Figures 2 and 4.

**9. On-line protocol, part II.** In this section we will briefly discuss the relation of our results about the IID model to Wilks's nonparametric prediction intervals and mention some relaxations of the on-line protocol.

*The univariate IID model.* The construction of prediction and tolerance intervals in the univariate IID model, which says that $y_1, y_2, \ldots$ form an IID sequence, was undertaken by many authors following the pioneering paper by Wilks (1941). Wilks's work was later extended to the multivariate case: see, for example, Fraser (1957); this extension, however, is not directly related to our IID predictors. For simplicity, let us assume in this subsection, as is customary in literature, that the distribution of one observation is continuous. Correspondingly, we will assume that the realized values of $y_n$, $n=1,2,\ldots$, are all different.

For each $n=1,2,\ldots$, define $T_n \in \{1,2,\ldots,n\}$ as the smallest $i$ such that $y_n < y_{(n-1,i)}$, where $y_{(n-1,1)}, \ldots, y_{(n-1,n-1)}$ is the sequence of the first $n-1$ observations $y_1, \ldots, y_{n-1}$ sorted in the ascending order; if $y_n > y_{(n-1,n-1)}$, set $T_n := n$. Each $T_n$ is a "pivot," being distributed uniformly on the set $\{1, \ldots, n\}$. Wilks suggested the following prediction intervals based on this fact: fix a number $r \in \{1, 2, \ldots\}$ and define $\Gamma_n^{2r/n}$, $n = 2r+1, 2r+2, \ldots$, to be the interval $(y_{(n-1,r)}, y_{(n-1,n-r)})$; the probability of error, $y_n \notin \Gamma_n^{2r/n}$, is then $2r/n$. Now Theorem 1 implies that the whole random sequence $(T_1, T_2, \ldots)$ has a known distribution: namely, it is distributed according to the product $U_1 \times U_2 \times \cdots$ of the uniform distributions $U_n$ on $\{1, \ldots, n\}$. In particular, Wilks' prediction intervals $\Gamma_n^{2r/n}$, $n = 2r+1, 2r+2, \ldots$, lead to independent errors.

*Relaxations of the on-line protocol.* This paper concentrates on the on-line prediction protocol. Smoothed conformal predictors lead to independent errors in the on-line protocol, and Theorem 2 suggests that conformal predictors are the most natural weakly valid confidence predictors. This is why we included the requirement of independence in the definition of strong validity, despite the fact that the error frequency can be shown to approach the error probability $\varepsilon$ with probability approaching one even when the requirement of independence is relaxed in certain ways.

The situation changes when we move outside the on-line protocol. The on-line protocol is natural, but in one respect it is overly restrictive: the true response $y_n$ becomes known before the prediction for the next response $y_{n+1}$ is made. It can be shown that the error frequency will still converge to $\varepsilon$ if the true response is only given for a small fraction of observations, and even for those observations it can be given with a delay [Vovk, Gammerman and Shafer (2005), Section 4.3; see also Vanderlooy, van der Maaten and Sprinkhuizen-Kuyper (2007) for a recent empirical study]. The independence of errors, however, will be lost (we can still have "approximate independence," but this is a much more elusive notion than ordinary independence).

**10. Conclusion.** In this paper we considered the problem of prediction in three main regression models. One of these models, the Gauss linear model, is the standard textbook one. The MVA model seems to have been somewhat neglected, partly because of philosophical reasons: according to the conditionality principle [Cox and Hinkley (1974), Section 2.3(iii)] one should condition on the observed values of the explanatory variables to make the prediction (or estimate, etc.) more relevant to the data at hand. In most of this paper we took a pragmatic approach, studying which models permit one to produce informative prediction intervals in different circumstances without being restricted a priori by general principles. We did use the sufficiency principle in our interpretation of Theorem 2, but we admit this makes the theorem less convincing. Surprisingly, the IID model appears to have been neglected in the field of regression, even in nonparametric statistics, where the value of this model is in principle well understood.

## APPENDIX: PROOFS OF THE THEOREMS

In this appendix we will prove the two main results stated in this paper, Theorems 1 and 2. A version of Theorem 1 was proved in Section 8.7 of Vovk, Gammerman and Shafer (2005), but we reproduce the principal points of the proof to make our exposition self-contained. A special case of Theorem 2 (namely, for the IID model) was proved in Section 2.6 of Vovk, Gammerman and Shafer (2005).

**Proof of Theorem 1.** In this proof, $\zeta_1, \zeta_2, \ldots$ will be random observations generated by $P \in \mathcal{P}$, $(\zeta_1, \zeta_2, \ldots) \sim P$, and $\tau_1, \tau_2, \ldots$ will be random numbers, $(\tau_1, \tau_2, \ldots) \sim U^\infty$. For each $n = 0, 1, \ldots$ let $\mathcal{G}_n$ be the $\sigma$-algebra generated by the random elements

$$S_n(\zeta_1, \ldots, \zeta_n), \zeta_{n+1}, \tau_{n+1}, \zeta_{n+2}, \tau_{n+2}, \ldots.$$

So $\mathcal{G}_0$ is the most informative $\sigma$-algebra and $\mathcal{G}_0 \supseteq \mathcal{G}_1 \supseteq \mathcal{G}_2 \supseteq \cdots$. It will be convenient to write $\mathbb{P}_\mathcal{G}(E)$ and $\mathbb{E}_\mathcal{G}(\xi)$ for the conditional probability $\mathbb{P}(E \mid \mathcal{G})$ and expectation $\mathbb{E}(\xi \mid \mathcal{G})$, respectively, given a $\sigma$-algebra $\mathcal{G}$.

LEMMA A.1. *For any step $n = 1, 2, \ldots$ and any $\varepsilon \in (0, 1)$,*

$$\mathbb{P}_{\mathcal{G}_n}(p_n \leq \varepsilon) = \varepsilon.$$

PROOF. For a given value of the summary $S_n(\zeta_1, \ldots, \zeta_n)$ of the first $n$ observations, consider the conditional distribution function $F$ of the random variable $\eta := A_n(S_{n-1}(\zeta_1, \ldots, \zeta_{n-1}), \zeta_n)$ (because of the total sufficiency, it does not matter whether we further condition on $\zeta_{n+1}, \tau_{n+1}, \zeta_{n+2}, \tau_{n+2}, \ldots$). Define $F(x-)$ to be $\sup_{t<x} F(t)$. Our task is to show that the conditional probability of the event

(A.1) $$1 - F(\eta) + \tau_n(F(\eta) - F(\eta-)) \leq \varepsilon$$

is $\varepsilon$ [since the left-hand side of (A.1) coincides with the right-hand side of the definition (2)]. The latter fact is usually stated in statistics textbooks for continuous $F$ [see, e.g., Cox and Hinkley (1974), page 66], but it is also easy to check in general. □

LEMMA A.2. *For any step $n = 1, 2, \ldots$, $p_n$ is $\mathcal{G}_{n-1}$-measurable.*

PROOF. This follows from the definition: $p_n$ is defined in terms of $\zeta_n$, $\tau_n$ and the summary of the first $n-1$ observations. □

Now we can easily prove the theorem. First we demonstrate that, for any $n = 1, 2, \ldots$ and any $\varepsilon_1, \ldots, \varepsilon_n \in (0, 1)$,

(A.2) $$\mathbb{P}_{\mathcal{G}_n}(p_n \leq \varepsilon_n, \ldots, p_1 \leq \varepsilon_1) = \varepsilon_n \cdots \varepsilon_1 \qquad \text{a.s.}$$

The proof is by induction on $n$. For $n = 1$, (A.2) is a special case of Lemma A.1. For $n > 1$ we obtain, from Lemmas A.1 and A.2, standard properties of conditional expectations, and the inductive assumption:

$$\begin{aligned}
\mathbb{P}_{\mathcal{G}_n}(p_n \leq \varepsilon_n, \ldots, p_1 \leq \varepsilon_1) &= \mathbb{E}_{\mathcal{G}_n}(\mathbb{E}_{\mathcal{G}_{n-1}}(\mathbb{I}_{p_n \leq \varepsilon_n} \mathbb{I}_{p_{n-1} \leq \varepsilon_{n-1}, \ldots, p_1 \leq \varepsilon_1})) \\
&= \mathbb{E}_{\mathcal{G}_n}(\mathbb{I}_{p_n \leq \varepsilon_n} \mathbb{E}_{\mathcal{G}_{n-1}}(\mathbb{I}_{p_{n-1} \leq \varepsilon_{n-1}, \ldots, p_1 \leq \varepsilon_1})) \\
&= \mathbb{E}_{\mathcal{G}_n}(\mathbb{I}_{p_n \leq \varepsilon_n} \varepsilon_{n-1} \cdots \varepsilon_1) \\
&= \varepsilon_n \varepsilon_{n-1} \cdots \varepsilon_1 \qquad \text{a.s.}
\end{aligned}$$

The "tower property" of conditional expectations immediately implies

$$\mathbb{P}(p_n \leq \varepsilon_n, \ldots, p_1 \leq \varepsilon_1) = \varepsilon_n \cdots \varepsilon_1.$$

Therefore, the distribution of the first $n$ p-values $p_1, \ldots, p_n$ is $U^n$, for all $n = 1, 2, \ldots$. This implies that the distribution of the infinite sequence $p_1, p_2, \ldots$ is $U^\infty$.

**Proof of Theorem 2.** In this proof, $Z := \mathbb{R}^K \times \mathbb{R}$ and $\zeta_i$ stands for $(\mathbf{x}_i, y_i)$. Let $n \in \mathcal{N}$.

For each summary $s \in \Sigma_n$ let $f(s)$ be the conditional probability given $S_n(\zeta_1, \ldots, \zeta_n) = s$ that $\Gamma$ makes an error at a significance level $\varepsilon$ when predicting $y_n$ from $\zeta_1, \ldots, \zeta_{n-1}$ and $\mathbf{x}_n$, the observations $\zeta_1, \zeta_2, \ldots$ being generated from $P \in \mathcal{P}$. We know that the expected value of $f(S_n(\zeta_1, \ldots, \zeta_n))$ is $\varepsilon$ under any $P \in \mathcal{P}$, and this, by the bounded completeness of $S_n$, implies that $f(s) = \varepsilon$ for almost all (under $PS_n^{-1}$ for any $P \in \mathcal{P}$) summaries $s$. Define $E(s, \varepsilon)$ to be the set of all pairs $(s', \zeta) = (s', (\mathbf{x}, y)) \in \Sigma_{n-1} \times Z$ such that $F_n(s', \zeta) = s$ (where $F_n$ is the function from the definition of the algebraic transitivity of the $S_n$) and $\Gamma$ makes an error at the significance level $\varepsilon$ when predicting $y$ and fed with $\zeta_1, \ldots, \zeta_{n-1}$ satisfying $S_{n-1}(\zeta_1, \ldots, \zeta_{n-1}) = s'$ and with $\mathbf{x}$ (since $\Gamma$ is invariant, whether an error is made depends only on $s'$, not on the particular $\zeta_1, \ldots, \zeta_{n-1}$). It is clear that

$$\varepsilon_1 \leq \varepsilon_2 \implies E(s, \varepsilon_1) \subseteq E(s, \varepsilon_2)$$

and

$$\mathbb{P}((S_{n-1}(\zeta_1, \ldots, \zeta_{n-1}), \zeta_n) \in E(s, \varepsilon) \mid S_n(\zeta_1, \ldots, \zeta_n) = s) = \varepsilon \qquad \text{a.s.,}$$

where $(\zeta_1, \zeta_2, \ldots) \sim P \in \mathcal{P}$.

In this proof we say "conformity measure" to mean a nonconformity measure which is used for computing p-values in the opposite way to (2): the ">" in (2) is replaced by "<." Let us check that the conformal predictor $\Gamma^\dagger$ determined by the conformity measure

$$A_n(s', \zeta) := \inf\{\varepsilon : (s', \zeta) \in E(F_n(s', \zeta), \varepsilon)\}$$

is at least as accurate as $\Gamma$. By the monotone convergence theorem for conditional expectations,

$$\mathbb{P}(A_n(S_{n-1}(\zeta_1, \ldots, \zeta_{n-1}), \zeta_n) \leq \varepsilon \mid S_n(\zeta_1, \ldots, \zeta_n) = s)$$
$$= \lim_{\delta \downarrow \varepsilon} \mathbb{P}(A_n(S_{n-1}(\zeta_1, \ldots, \zeta_{n-1}), \zeta_n) < \delta \mid S_n(\zeta_1, \ldots, \zeta_n) = s)$$
$$\leq \lim_{\delta \downarrow \varepsilon} \mathbb{P}((S_{n-1}(\zeta_1, \ldots, \zeta_{n-1}), \zeta_n) \in E(s, \delta) \mid S_n(\zeta_1, \ldots, \zeta_n) = s)$$
$$= \lim_{\delta \downarrow \varepsilon} \delta = \varepsilon \qquad \text{a.s.,}$$

where $(\zeta_1, \zeta_2, \ldots) \sim P \in \mathcal{P}$ and $\delta$ is constrained to be a rational number. Therefore, at each significance level $\varepsilon$ and for all $(\zeta_1, \ldots, \zeta_n) \in Z^n$,

$$\begin{aligned}
y_n &\in (\Gamma^\dagger)^\varepsilon(\zeta_1, \ldots, \zeta_{n-1}, \mathbf{x}_n) \\
&\iff \mathbb{P}(A_n^{\text{rnd}} \leq A_n^{\text{obs}} \mid S_n^{\text{rnd}} = S_n^{\text{obs}}) > \varepsilon \\
&\implies A_n^{\text{obs}} > \varepsilon \\
&\implies (S_{n-1}(\zeta_1, \ldots, \zeta_{n-1}), \zeta_n) \notin E(S_n(\zeta_1, \ldots, \zeta_n), \varepsilon) \\
&\iff y_n \in \Gamma^\varepsilon(\zeta_1, \ldots, \zeta_{n-1}, \mathbf{x}_n) \qquad \text{a.s.,}
\end{aligned}$$

in the notation of (2) and for $(\xi_1, \xi_2, \ldots) \sim P \in \mathcal{P}$.

**Acknowledgments.** We have greatly benefited from Glenn Shafer's advice and from a discussion with Steffen Lauritzen. Comments by the anonymous referees and Professor Susan Murphy helped us improve the presentation.

**Note added in proof.** The R package PredictiveRegression, available from CRAN, implements the three prediction algorithms (IID predictor, Gauss predictor and MVA predictor) described in this paper.

V. Vovk
I. Nouretdinov
A. Gammerman
Computer Learning Research Centre
Department of Computer Science
Royal Holloway, University of London
Egham, Surrey TW20 0EX
United Kingdom
E-mail: vovk@cs.rhul.ac.uk